\documentclass[11pt]{article}

\usepackage{amssymb} \usepackage{latexsym} \usepackage{amsthm}
\usepackage{exscale}  \usepackage{amsmath}
\usepackage{epsf}

\newtheorem{theorem}{Theorem}[section]
\newtheorem{lemma}[theorem]{Lemma}
\newtheorem{proposition}[theorem]{Proposition}

\newenvironment{remark}{\stepcounter{remark}\noindent $\triangleleft$
 \textbf{Remark \thesection.\theremark: }}{$\triangleright$\newline}
\newcounter{remark}[section]

\newcommand{\eps}{\varepsilon}
\newcommand{\nquad}{\negthickspace\negthickspace
\negthickspace\negthickspace}

\title{BV instability for the Lax-Friedrichs scheme}
\author{Paolo Baiti \thanks{Dipartimento di Matematica e Informatica,
	Universit\`a di Udine, Via delle Scienze 206,
	Udine 33100, Italy; E-mail: \texttt{baiti@dimi.uniud.it}}
	\and
	Alberto Bressan\thanks{Department of Mathematics, Penn State
	University,
	UP, PA 16802, U.S.A.;
	E-mail: \texttt{bressan@math.psu.edu}}
	\and
	Helge Kristian Jenssen\thanks{Department of Mathematics, North
	Carolina State University,
	Raleigh, NC 27695, U.S.A.;
	E-mail: \texttt{hkjensse@math.ncsu.edu}. Research supported in part
	by the NSF under grant DMS-0206631}
	}

\date{}

\begin{document}

\maketitle

\abstract{It is proved that discrete shock profiles (DSPs) for the
Lax-Friedrichs scheme for a system of conservation laws do not necessarily
depend continuously in BV on their speed. We construct examples
of $2 \times 2$-systems for which there are sequences of DSPs with speeds converging
to a rational number. Due to a resonance phenomenon, the difference between
the limiting DSP and any DSP in the sequence will contain an order-one
amount of variation.}

\section{Introduction}
Consider a strictly hyperbolic $n\times n$ system of conservation laws
in one space dimension:
\begin{equation}\label{1.1}
	u_t+f(u)_x=0.
\end{equation}
For initial data with small total variation,
the existence of a unique entropy weak solution
is well known \cite{Glimm65}, \cite{BressanLiuYang99}, \cite{BressanBook}.
A closely related question is the stability and convergence of various types
of approximate solutions.
For vanishing viscosity approximations
\begin{equation}\label{1.2}
	u_t+f(u)_x=\eps\, u_{xx}\,,
\end{equation}
uniform BV bounds, stability and convergence as $\eps\to 0$
were recently established in \cite{BianchiniBressan04}.
Assuming that all the eigenvalues of the Jacobian matrix $Df(u)$
are strictly positive,
similar results are also proved in \cite{Bianchini03}
for solutions constructed by the semidiscrete (upwind)
Godunov scheme
$$
\frac{d}{dt}\, u_j(t)+ \frac{1}{\Delta x}\Big[
f\big(u_j(t)\big)-f\big(u_{j-1}(t)\big)\Big]=0\,,
\qquad u_j(t)=u(t,j\,\Delta x)\,.
$$
In the present paper we study the case of fully discrete
schemes, where the derivatives w.r.t.~both 
time and space are replaced by
finite differences. 

We recall that, for the $2\times 2$ system of isentropic gas dynamics,
the convergence of Lax-Friedrichs and Godunov approximations was 
proved in \cite{DingChenLuo89}, within the framework of compensated compactness.
Further results have been obtained for {\sl straight line systems},
where all the Rankine-Hugoniot 
curves are straight lines.  For straight line systems of two equations LeVeque and Temple 
\cite{LeVequeTemple85}
utilized the existence of Riemann invariants to prove stability and convergence 
of the Godunov scheme.
For $n\times n$-systems in the same class, uniform BV bounds, stability and convergence
of a relaxation scheme, Godunov and Lax-Friedrichs approximations were established in
\cite{BressanShen00}, \cite{BressanJenssen00}, \cite{YangZhaoZhu03}, respectively.  
The analysis relies on the fact that, due to the very
particular geometry, the interaction of waves of the same family
does not generate additional oscillations.

A key ingredient in the arguments in \cite{BianchiniBressan04} and \cite{Bianchini03} is the local decomposition
of the approximate solutions in terms of travelling waves.
To achieve a good control on the new waves produced by  
interactions of waves of a same family, it is essential that
the center manifold of travelling profiles has a certain degree of smoothness.
We show in this paper that this smoothness is lacking in the case of fully discrete schemes.
As remarked by Serre \cite{Serre96}, for general hyperbolic systems 
the discrete shock profiles cannot depend continuously on the speed $\sigma$,
in the BV norm.  In the present paper we construct an explicit example 
showing how this happens.

Our basic example is provided by a 
$2\times 2$ system in triangular form
\begin{align}
     u_t + f(u)_x &= 0, \label{eq1}\\
     v_t + g(u)_x &= 0.\label{eq2}
\end{align}
The characteristic speeds are $0$ and $f'(u)$ and the system is
strictly hyperbolic provided $f'(u)>0$. 
The Lax-Friedrichs scheme for (\ref{eq1})-(\ref{eq2}) with $\Delta
x=\Delta t$ takes the form
\begin{align}
  u_{n+1,j} &= \frac{1}{2}\big(u_{n,j+1}+u_{n,j-1}\big) -
  \frac{1}{2}\big(f(u_{n,j+1})-f(u_{n,j-1})\big),  \label{LF1}\\
  v_{n+1,j} &= \frac{1}{2}\big(v_{n,j+1}+v_{n,j-1}\big) -
  \frac{1}{2}\big(g(u_{n,j+1})-g(u_{n,j-1})\big).  \label{LF2}
\end{align}
For the rest of the paper we fix a flux  
function $f(u)$ which satisfies $f'(u)>1/4$, say, and also the CFL condition
$|f'(u)|<1$, for all $u\in\mathbb R$. We also take the other flux 
function $g$ constant outside a bounded interval.

A {\em discrete shock profile} (DSP) with speed $\lambda$ for
(\ref{LF1})-(\ref{LF2}) is a pair of functions
\[\big(U(x),V(x)\big)=\big(U^{(\lambda)}(x),V^{(\lambda)}(x)\!\big)\] 
satisfying 
\begin{align}
  U(x-\lambda) &= \frac{U(x+1)+U(x-1)}{2} - 
  \frac{f(U(x+1)\!)-f(U(x-1)\!)}{2}, \label{DSP1}\\ \nonumber \\
  V(x-\lambda) &= \frac{V(x+1)+V(x-1)}{2} - 
  \frac{g(U(x+1)\!)-g(U(x-1)\!)}{2}.  \label{DSP2}
\end{align}
We will give a rational speed $\lambda=p/q$ and a sequence of rational
perturbations $\eps_n$ for which the difference
$V^{(\lambda)}-V^{(\lambda+\eps_n)}$ contains an $O(1)$
amount of variation, uniformly with respect to $n$. The variation occurs far downstream in
 an interval of the form
$[-C\eps_n^{-2},-c\eps_n^{-2}]$.
%
%
\section{Outline of construction} 
As a motivation for the later computations we consider the easier case of the heat equation
with a point-source. If the source acts continuously in time
and concentrated along the line $x=\sigma t$, then
\begin{equation}\label{2.4}
	v_t-v_{xx}=\delta_{t,\sigma t}
\end{equation}
with $\sigma>0$. In this case one finds the travelling wave solution
\begin{equation}\label{2.5}
	v(t,x)=\phi(x-\sigma t)
\end{equation}
with
\begin{equation}\label{2.6}
	\phi(y)=\int_0^\infty G(t,y+\sigma t)\,dt =\left\{
	\begin{array}{ll}
	\sigma^{-1}e^{-\sigma y}\quad & \mbox{if $y\geq 0$},\\
	\sigma^{-1}\quad & \mbox{if $y\leq 0$.}\end{array}\right.
\end{equation}
Here $G(t,x):= e^{-x^2/4t}/2\sqrt{\pi t}$ is the standard heat kernel.
Notice that the travelling profile can also be obtained as
the value at time $t=0$ of a solution of (\ref{2.4})
defined for $t\in \,]-\infty,\,0]$. We also have
\begin{equation}\label{2.7}
	\phi'(y)=\int_0^\infty G_x(t,y+\sigma t)\,dt =\left\{
	\begin{array}{ll}
	-e^{-\sigma y}\quad &\mbox{if $y>0$,}\\
	0\quad & \mbox{if $y< 0$.}\end{array}\right.
\end{equation}
%
\begin{figure}[ht]\label{fig1}
\centerline{\epsfxsize8truecm\epsfbox{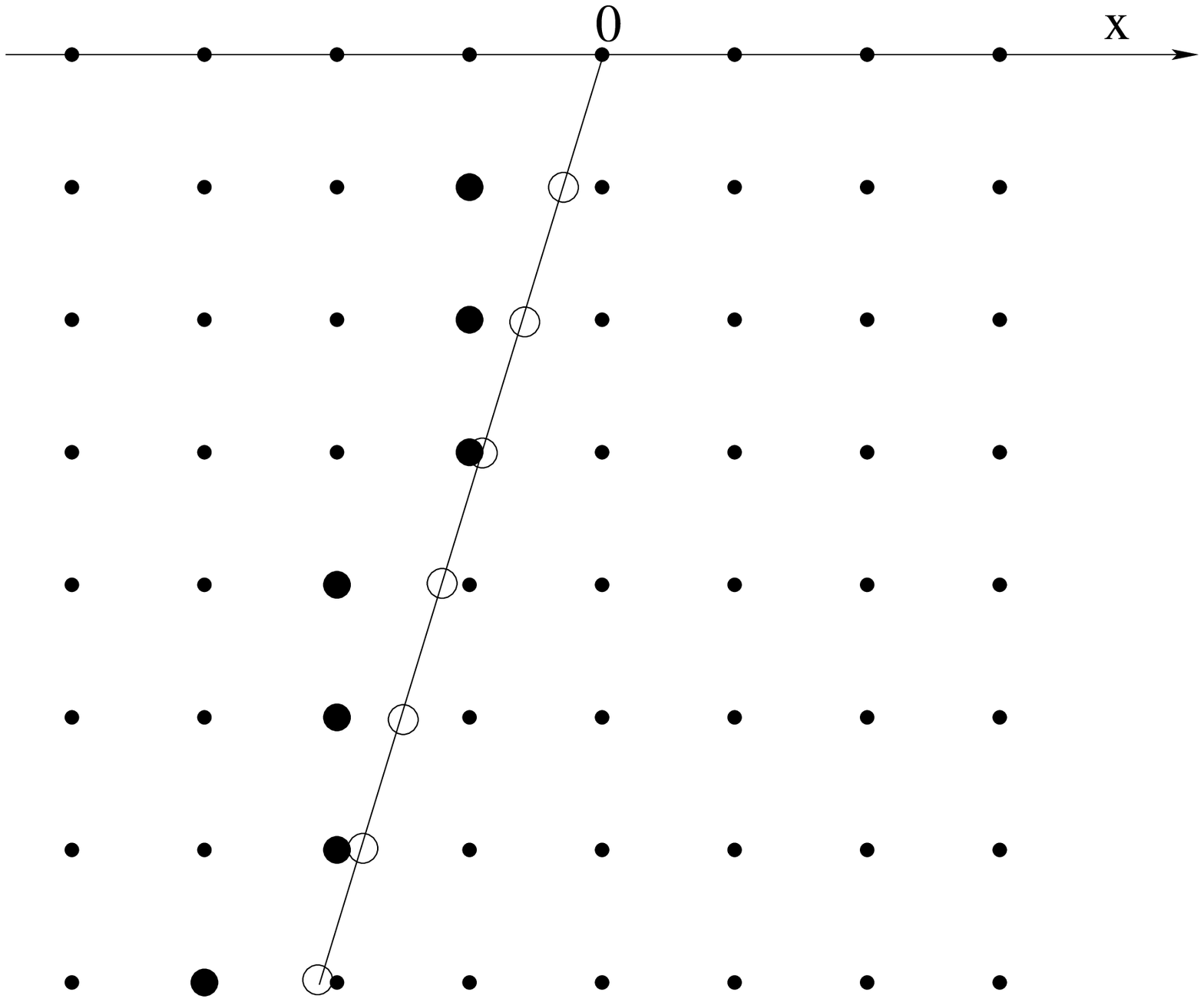}}
\centerline{Figure~1}
\end{figure}
%
Next, we consider the case 
where the sources are located on a discrete set of points 
$P_n=(n,\sigma n)$,
with $n$ integer (the white circles in Figure~1)
\begin{equation}\label{2.8}
	v_t-v_{xx}=\delta_{n,\sigma n}\,.
\end{equation}
We again assume that $\sigma>0$ and  consider a solution of (\ref{2.8})
defined for $t\in ]-\infty,\,0]$.
Its value at time $t=0$ is now computed as
\begin{equation}\label{2.9}
	v(0,y)=\Phi(y)\doteq\sum_{n\geq 1} G(n, y+\sigma n)\,.
\end{equation}
For $y\to \infty$, it is clear that $\Phi(y)$ tends to zero
exponentially fast, together with all derivatives.
We wish to understand how the oscillations decay for $y\to -\infty$, 
i.e.~far downstream from the shock. 
For $y<0$, the sum (\ref{2.9}) can be expressed as an integral
\begin{equation}\label{2.10}
	\Phi(y)\doteq\sum_{n\geq 1} G(n,\, y+\sigma n)=
	\int_0^\infty G(t,y+\sigma t)\, \big(1+h_1'(t)\big)\,dt,
\end{equation}
where
$$h_1(t)\doteq [\![t ]\!]-t+1/2\,.$$
By induction, we can find a sequence of periodic functions
$h_m$ such that (Figure~2)
$$h_m(t)=h_m(t+1)\,,\qquad\int_0^1 h_m(t)\,dt=0\,,\qquad
\frac{d}{dt}h_m(t)=h_{m-1}(t)\,.$$
 %
\begin{figure}[ht]
\centerline{\epsfxsize8truecm\epsfbox{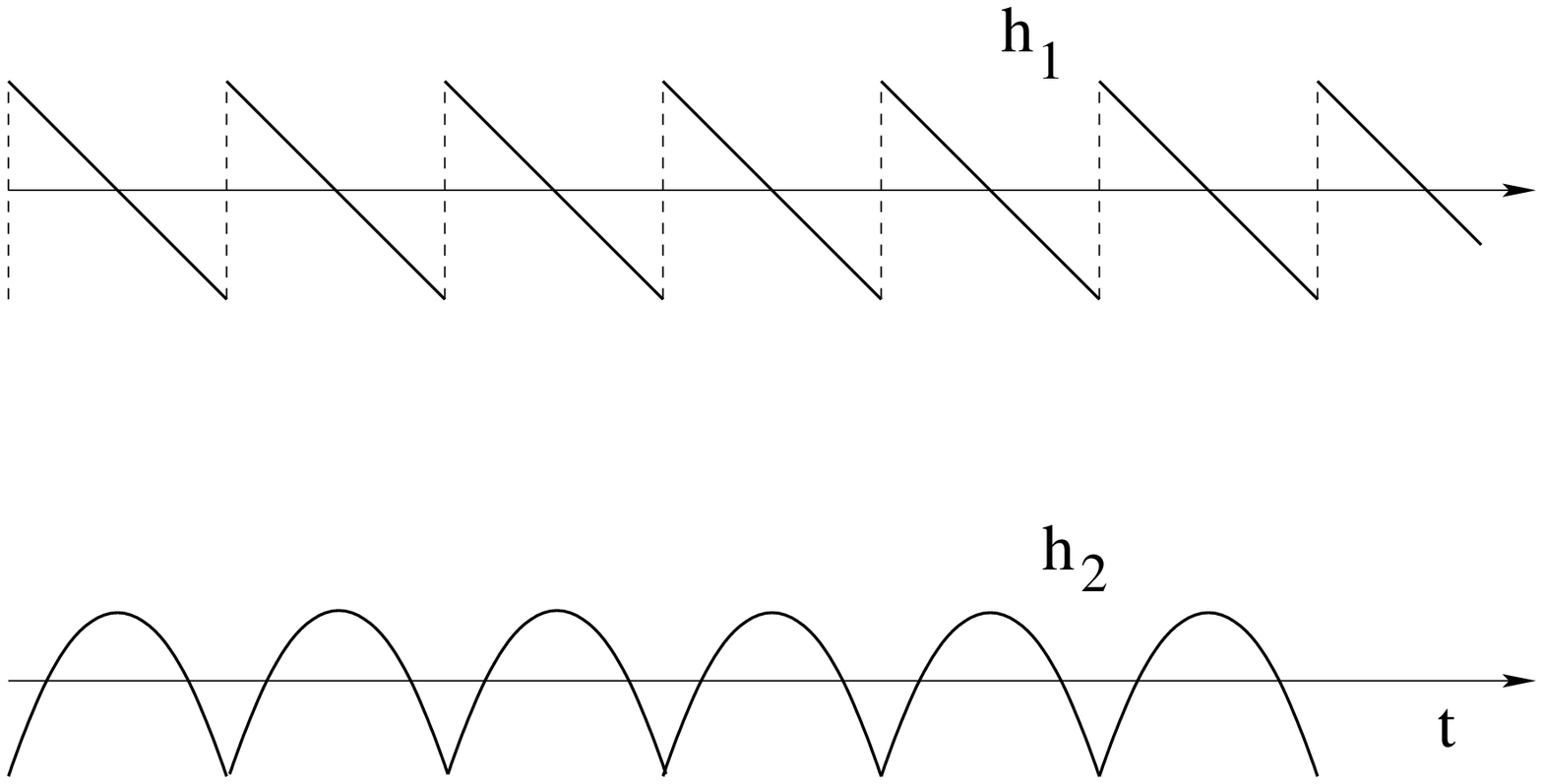}}
\centerline{Figure~2}
\end{figure}
Integrating by parts and recalling (\ref{2.6}), from (\ref{2.10}) we obtain
\begin{align}
	\Phi(y)&= \int_0^\infty G(t,y+\sigma t)\,\Big(1+ \frac{d^m}{dt^m}
	h_m(t)\Big)\,dt\\	
	&= \frac{1}{\sigma}+(-1)^m\int_0^\infty \frac{d^m}{dt^m}
	G(t,y+\sigma t) \,h_m(t)\,dt\,.\label{2.11}
\end{align}
Since
$$h_m(t)=\int_{\xi_m}^t h_{m-1}(s)\,ds,\qquad\qquad t\in [0,1],$$
for some point $\xi_m\in [0,1]$, by induction we find
$$\big|h_m(t)\big|\leq 1
\qquad\qquad\forall\, m\geq 1,\,\forall\, t\in\mathbb R\,.$$
The identities
$$G(t,x)=t^{-1/2}G(1,x/\sqrt t)\,,\qquad\qquad G_t=G_{xx}\,,$$
imply
\begin{align}
	\frac{\partial^m}{\partial x^m}G(t,x)&=t^{-(m+1)/2}
	\cdot\frac{\partial^m}{\partial x^m}G(1,x/\sqrt t)\,,\\
	\frac{\partial^m}{\partial t^m}G(t,x)&=t^{-(2m+1)/2}
	\cdot\frac{\partial^m}{\partial t^m}G(1,x/\sqrt t)\,.
\end{align}
We recall here some basic estimates for the heat kernel and its derivatives, which
we will use throughout the paper. For $k=0,1,\dots$ we have
\begin{align}
\left\|\frac{\partial^k}{\partial x^k}G(t,\cdot)\right\|_{\mathbf{L}^{\infty}}
&=O(1)\cdot\frac{1}{t^{(k+1)/2}}, \label{Gx_esti}\\
\left\|\frac{\partial^k}{\partial t^k}G(t,\cdot)\right\|_{\mathbf{L}^{\infty}}
&=O(1)\cdot\frac{1}{t^{(2k+1)/2}}. \label{Gt_esti}
\end{align}
In addition we observe that, as $y\to -\infty$, the function
$t\mapsto G(t,\,y+\sigma t)$ becomes exponentially small
together with all its derivatives, outside the interval centered at
$|y|/\sigma$ with width $|y|^{\delta+1/2}\,$, for any $\delta>0$.  
More precisely
\begin{equation}\label{2.12}
	\sup_{|t+y/\sigma|<|y|^{\delta+1/2}}
	\left|\frac{d^m}{ dt^m}
	G(t,y+\sigma t) \right|=O(1)\cdot e^{c_\delta y}\qquad
	\hbox{as }y\to -\infty\,,
\end{equation}
for some constant $c_\delta>0$.

Letting $y\to-\infty$, for every $m\geq 1$ the above estimates imply
\begin{equation}\label{2.13}
	\left|\Phi(y)-\frac{1}{\sigma}\right| 	\leq 
	\int_0^\infty \left|\frac{d^m}{dt^m}
	G(t,y+\sigma t) \right|\,dt=O(1)\cdot y^{-m/2}.
\end{equation}
Similarly,
\begin{equation}\label{2.14}
	\big|\Phi'(y)\big|\leq 
	\int_0^\infty \left|\frac{d^{m+1}}{dt^{m+1}}
	G(t,y+\sigma t) \right|\,dt=O(1)\cdot y^{-(m+1)/2}.
\end{equation}
Since $m\geq 1$ is arbitrary, this shows that
the function $\Phi'$ is rapidly decreasing as $y\to -\infty$.
In particular, taking $m=2$ in (\ref{2.14}) we obtain the integrability of 
$\Phi'$, hence a bound on the total variation of $\Phi$.

Finally, assume that the impulses are located not at the points
$P_n=(n,\sigma n)$ but at the points with integer coordinates
$Q_n\doteq \big(n,[\![\sigma n]\!] \big)$ (the black circles in Figure~1)
\begin{equation}\label{2.15}
	v_t-v_{xx}=\delta_{n,[\![\sigma n]\!]}\,,
\end{equation}
where $[\![a]\!]$ denotes the integer part of a real number $a$.

Again we consider a solution defined for $t\in \,]-\infty,\,0]$
and study its profile at the terminal time $t=0$. Assuming $\sigma>0$, a direct computation yields
$$v(0,y-1)=\Psi(y):=
\sum_{n\geq 1} G\big(n,y+[\![\sigma n]\!]\big).$$
Because of (\ref{2.13}), to determine the asymptotic behavior as $y\to-\infty$,
it suffices to estimate the difference
$$
K(y)\doteq\Psi(y)-\Phi(y)=
-\sum_{n\geq 1} \Big[G(n,y+\sigma n\big)
-G\big(n,y+[\![\sigma n]\!]\big)\Big]\,.
$$
It is here that, if the speed $\sigma$ is close to a rational,
a resonance is observed.   To see a simple case, let
$\sigma=1+\eps$, with $\eps>0$ small.  Then
we can approximate
\begin{align}
	K(y)&\approx -\sum_{n\geq 1} G_x(n,y+\sigma n\big)
	\big(\sigma n-[\![\sigma n]\!]\big) \notag\\
	&\approx -\int_0^\infty G_x(t,y+\sigma t)\, \big(\eps t -[\![\eps t]\!]\big)\,dt\,.\label{2.16}
\end{align}
The functions appearing in the above integration are shown in Figure~3.
We recall that 
$$
\int_0^\infty G_x(t,y+\sigma t)\,dt=
-\int_0^\infty \frac{y+\sigma t}{4t\sqrt{\pi t}}\, 
\exp\left\{ -\frac{(y+\sigma t)^2}{4t}\right\}\,dt=0
$$
for every $y<0$.
Set $y_\eps\doteq -\eps^{-2}$.
When $y$ ranges within the interval
$$
I_\eps\doteq [y_\eps,y_\eps/2]=[-\eps^{-2},-\eps^{-2}/2]\,,
$$
the integral in (\ref{2.16}) can be of the same order of magnitude as
$$
\int_0^\infty \big|G_x(t,y_\eps+\sigma t)\big|\,dt\geq
 c_0\, y_\eps^{-1/2}=c_0\,\eps\,.
 $$
%
\begin{figure}[ht]\label{fig3}
\centerline{\epsfxsize8truecm\epsfbox{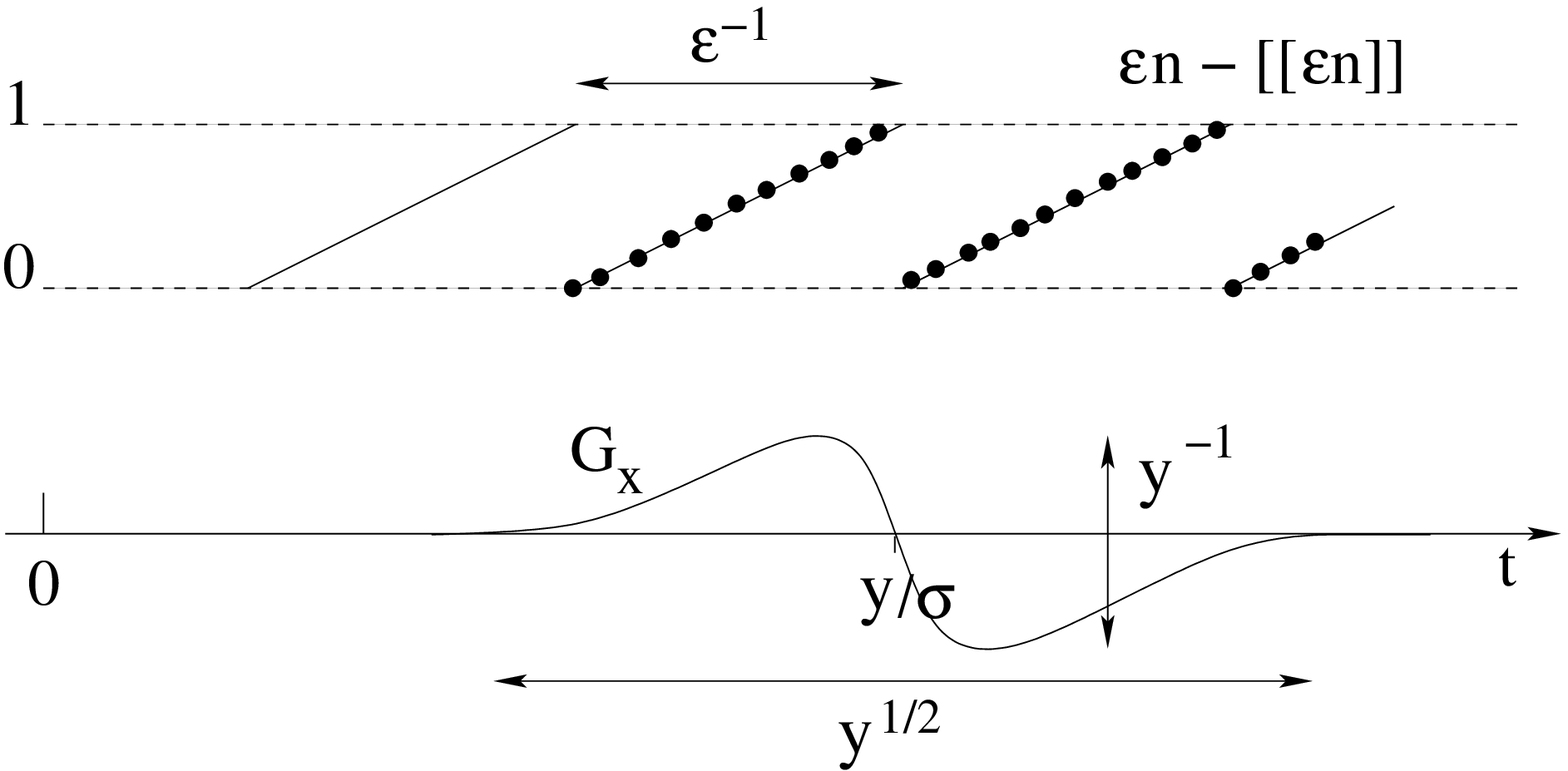}}
\centerline{Figure~3}
\end{figure}
Moreover, each time that
$y$ increases by an amount $\Delta y=\eps^{-1}$, the phase of the
fractional part $[[\eps y]]-\eps y$ goes through a full cycle, hence the map
$$y\mapsto \int_0^\infty G_x(t,y+\sigma t)\, \big(
\eps t -[[\eps t]]\big)\,dt$$
oscillates by an amount $\geq c_1\eps$.  In all, 
we have approximately
$1/2\eps$ cycles
within the interval $I_\eps$. Hence the total variation of the 
discrete profile
$\Psi^{(1+\eps)}$ on $I_\eps$ can be estimated as
\begin{equation}\label{2.17}
	T.V.\big\{\Psi^{(1+\eps)}~;~I_\eps \big\}\geq c_2 
\end{equation}
for some constant $c_2>0$ independent of $\eps$.
We write here $\Psi=\Psi^{(1+\eps)}$ to 
emphasize that the profile depends on the speed $\sigma=1+\eps$.
By (\ref{2.17}) it is clear that, as $\eps\to 0+$, the 
functions $\Psi^{(1+\eps)}$
do not form a Cauchy sequence and cannot converge in the space BV.
%
%
\section{Discrete shock profiles for the system}
Since (\ref{eq1}) is a decoupled scalar equation, the Lax-Friedrichs scheme 
(\ref{LF1}) for this 
component admits discrete shock profiles connecting any two states, see Jennings \cite{Jennings74}. 
For given left and right states 
$u_-,u_+$ we denote by $U^\lambda(x)$ the DSP connecting $u_-$ to $u_+$
and moving with speed 
\[\lambda=\frac{[f(u)]}{[u]}.\]
We proceed to construct the second component of a DSP with speed $\lambda$ by
using the discrete Green kernel for (\ref{LF2}) and the Duhamel principle.
Consider the Lax-Friedrichs scheme for the equation $z_t=0$,
\begin{equation}\label{linear_LF}
	z_{n+1,j} = \frac{z_{n,j+1}+z_{n,j-1}}{2}.
\end{equation}
We observe that the discrete Green's function $K_{n,k}$ for (\ref{linear_LF}) is given by
\begin{equation}\label{K}
	K_{n,k}=\left\{\begin{array}{ll}
		\left(\frac{1}{2}\right)^n \binom{n}{(n-k)/2} &
	\textrm{for $k=-n,-n+2,\dots,n-2,n$},\\
	 0 &\textrm{otherwise}.\end{array}\right.
\end{equation}
Given any DSP for the first equation, a DSP for (\ref{LF2}) is obtained by 
prescribing vanishing $v$-data at time
$n=-\infty$, and then letting the $u$-terms act as a source in (\ref{LF2}) 
from $n=-\infty$ to $n=0$. Consider first the difference equation
\begin{equation}\label{2nd_comp}
   v_{n+1,j} = \frac{1}{2}\big(v_{n,j+1}+v_{n,j-1}\big)+\psi_{n,j}, 
\end{equation}
where the sources $\psi_{n,j}$ are assumed given.
By Duhamel's principle we have that if vanishing data are given
at time step $-N$, then the solution of (\ref{2nd_comp}) at time step
$n\geq -N$ is  
\begin{equation}\label{solution}
     v_{n,j} =   \sum_{m=-N-1}^{n-1} \sum_{k\in\mathbb Z} \psi_{m,j-k}K_{n-1-m,k}.
\end{equation}
To apply this in our situation we introduce the functions
$\psi,\ H:\mathbb R\to \mathbb R$ by 
\[\psi(s):=\frac{d}{ds} g(U(s)\!),\]
and
\[H(x):=-\frac{1}{2}\big[g(U(x+1)\!)-g(U(x-1)\!)\big]=
-\frac{1}{2}\int_{x-1}^{x+1}\psi(s)\, ds ,\]  
where $U=U^{(\lambda)}$ is a scalar DSP for the first equation. In this
case (\ref{LF2}) may be written in the form (\ref{2nd_comp}) with the
sources $\psi_{n,j}$ given by 
\[\psi_{n,j}=H(j-\lambda n).\]
From now on we make the assumption that $g(u)$ is such that $\psi$ and 
hence also $H$ have compact support.
\begin{proposition}\label{2nd_DSP}
  	The pair of functions $(U^{(\lambda)}, V^{(\lambda)})$
	where $U^{(\lambda)}$ is a DSP for (\ref{LF1}) and
  	$V^{(\lambda)}$ is defined by  
  	\[V^{(\lambda)}(x):=\sum_{n=1}^{\infty}\sum_{k\in\mathbb
  	Z}H(x-k+\lambda n)K_{n-1,k},\]
  	is a DSP for the system (\ref{LF1})-(\ref{LF2}).
\end{proposition}
\begin{remark} The proof of this proposition is immediate once it is verified that
the double sum converges. Since $H$ has compact 
support the convergence follows.
\end{remark}
We next give a useful integral representation of $V^{(\lambda)}(x)$. For a
fixed $\xi$ we define the function $v^{(\lambda)}(\cdot;\xi):\mathbb 
R\to \mathbb R$ by  
\[v^{(\lambda)}(x;\xi):=\sum_{n=1}^{\infty}\left(K_{n-1,[\![x+\lambda n
    -\xi]\!]}+K_{n-1,[\![x+\lambda n -\xi]\!]+1}\right).\]
\begin{proposition}\label{DSP_repr}
	The function $V^{(\lambda)}(x)$ is given by
	\begin{equation}\label{repr}
	V^{(\lambda)}(x)=-\frac{1}{2}\int_{-\infty}^{\infty}\psi(\xi)
	v^{(\lambda)}(x;\xi)\, d\xi.
	\end{equation}
\end{proposition}
\begin{proof}
  \begin{align*}
    &V^{(\lambda)}(x) = \sum_{n=1}^{\infty}\sum_{k\in\mathbb Z} H(x-k+\lambda
    n)K_{n-1,k} \\ 
    &= -\frac{1}{2} \sum_{n=1}^{\infty}\sum_{k\in\mathbb
    Z}\left[ \int_{x-k+\lambda n}^{x-k+\lambda n+1}\psi(\xi) \, d\xi
    K_{n-1,k} +
    \int_{x-k+\lambda n-1}^{x-k+\lambda n} \psi(\xi) \, d\xi K_{n-1,k}
    \right]\\
    &= -\frac{1}{2}\sum_{n=1}^{\infty}\sum_{k\in\mathbb Z} \left[
    \int_{x-k+\lambda n}^{x-k+\lambda n+1}\psi(\xi)K_{n-1,[\![x+\lambda
    n -\xi]\!]+1} \, d\xi \right.\\
    &\qquad\qquad\qquad\qquad\qquad\qquad\qquad
    \left.+\int_{x-k+\lambda n-1}^{x-k+\lambda n} \psi(\xi) 
    K_{n-1,[\![x+\lambda n -\xi]\!]} \, d\xi
    \right]\\
    &= -\frac{1}{2}\sum_{n=1}^{\infty}\left[\int_{-\infty}^{\infty}
    \psi(\xi)K_{n-1,[\![x+\lambda n -\xi]\!]+1} \, d\xi  +
    \int_{-\infty}^{\infty} \psi(\xi) 
    K_{n-1,[\![x+\lambda n -\xi]\!]} \, d\xi\right]\\
    &= -\frac{1}{2}\int_{-\infty}^{\infty}
    \psi(\xi)\sum_{n=1}^{\infty} \left(K_{n-1,[\![x+\lambda n -\xi]\!]+1}  
    + K_{n-1,[\![x+\lambda n -\xi]\!]} \right)\, d\xi\\
    &= -\frac{1}{2}\int_{-\infty}^{\infty}\psi(\xi)v^{(\lambda)}(x;\xi)\, d\xi.
  \end{align*}
\end{proof}
%
%
\section{Approximation of the DSP in terms of the heat
  kernel} 
We will compare solutions of the Lax-Friedrichs scheme with
certain solutions of the heat equation. As a first step we approximate
the discrete Green's function $K_{n,k}$ using the heat kernel 
\[G(t,x) = \frac{1}{2\sqrt{\pi t}}e^{-\frac{x^2}{4t}}.\]
We use
the following notation (see \cite{Feller68}):  
\[a_k(\nu):=\left(\frac{1}{2}\right)^{2\nu}\binom{2\nu}{\nu+k}.\]
By Stirling's formula we have
\[a_k(\nu)=h\mathcal N(hk)\cdot \exp(\eps_1-\eps_2),\]
where 
\[h=\sqrt{\frac{2}{\nu}},\qquad \mathcal N(x)=\frac{1}{\sqrt{2\pi}}
e^{-\frac{x^2}{2}},\] 
and the errors $\eps_1,\eps_2$ satisfy
\[-\frac{3k^2}{4\nu^2}<\eps_1<\frac{k^4}{4\nu^3},\qquad
\textrm{provided $|k|<\nu/3$, and }\eps_2=O(1/\nu).\]  
There are two cases to consider depending on whether both $n$ and $k$
are even or both are odd. 
\begin{itemize}
  \item{} Case 1: $n=2m$, $k=2l$. In this case we have, 
    \[K_{n,k}=a_l(m)=2G\left(\frac{n}{2},k\right)
    \cdot e^{\eps_1-\eps_2}.\]  
  \item{} Case 2: $n=2m+1$, $k=2l+1$. In this case we have,
    \begin{align*}
      K_{n,k}&=\frac{n}{n+k}a_l(m)\\
      &=2\frac{n}{n+k}G\left(\frac{n-1}{2},k-1\right) 
      \cdot e^{\eps_1-\eps_2}\\
      &=2G\left(\frac{n}{2},k\right)\cdot e^{\eps_1-\eps_2}
      \frac{n}{n+k}\sqrt{\frac{n}{n-1}}
      e^{\frac{k^2}{2n}-\frac{(k-1)^2}{2(n-1)}}.
    \end{align*}
\end{itemize}
We will use this approximation only in the case when   
$|k|\leq O(1)n^{1/2+\delta}$ where $0<\delta\ll 1$. It follows that in either case we have
\[e^{\eps_1-\eps_2}=1+O(1)n^{-1+4\delta}.\]
An easy calculation shows that under the same condition on $k$,
\[\frac{n}{n+k}\sqrt{\frac{n}{n-1}}
e^{\frac{k^2}{2n}-\frac{(k-1)^2}{2(n-1)}} =1+O(1)n^{-1+\delta}.\]  
Recalling (\ref{Gx_esti}), summing up we have the following.
\begin{proposition}
For $n\geq 1$ and for $|k|\leq O(1)n^{1/2+\delta}$, with $0<\delta\ll 1$, we have
\begin{equation}\label{approx}
  	K_{n,k}=
  	2G\left(\frac{n}{2},k\right)+O(1)n^{-3/2+\delta}. 
\end{equation}
\end{proposition}
For a given speed $\lambda>0$ and for $\delta\in]0,1/2[$ 
we define the time interval
\begin{equation} \label{interval}
	I(y;\lambda,\delta)=\left[\frac{y}{\lambda}-y^{1/2+\delta},\frac{y}{\lambda}+y^{1/2+\delta}\right].
\end{equation}
We will make repeated use of the fact that, for $y<0$, the indices outside 
$I(|y|;\lambda,\delta)$ contribute 
exponentially little to the sum
\[S(y):=\sum_{n=1}^{\infty}K_{n,[\![y+\lambda n]\!]}.\]
\begin{proposition}\label{tails}
For  $y<0$, $\lambda>0$ and $\delta\in]0,1/2[$ we have
\begin{equation}\label{tail}
	\sum_{n\notin I(|y|;\lambda,\delta)}K_{n,[\![y+\lambda n]\!]} \leq
	O(1)e^{-C(\delta,\lambda) |y|^{c(\lambda,\delta)}},  
\end{equation}
for some positive constants $C(\delta,\lambda)$ and $c(\lambda,\delta)$.
\end{proposition}
\begin{proof}
For notational convenience assume that $\lambda=1$, the case $\lambda<1$
being similar. First of all we have 
$K_{n,[[y+\lambda n]]}\neq 0$ iff $n\geq \frac{|y|}{2}$.
We divide the above sum in four parts, where $n$ ranges over
$I_{1}=\big]\frac{|y|}{2},\frac{|y|+1}{2}\big]$,
$I_{2}=\big]\frac{|y|+1}{2},|y|-|y|^{1/2+\delta}\big[$,
$I_{3}=\big]|y|+|y|^{1/2+\delta},2|y|\big[$ and
$I_{4}=\big[2|y|,+\infty\big[$, respectively.
If $n\in I_{1}$ (and there is 
at most one such $n$) then $[\![y+n]\!]=-n$ and $n=O(|y|)$, so that
$K_{n,[[y+ n]]}=\frac{1}{2^{n}}=2^{-C |y|}$ is transcendentally small.
For the remaining indexes we can use the following estimate obtained by 
Stirling's formula
\begin{equation}\label{stirling}
	\binom{n}{k}\leq \frac{Cn^{n+1/2}}{k^{k+1/2}
	(n-k)^{n-k+1/2}},\qquad\textnormal{for }0<k<n.
\end{equation} 
If $n\in I_{2}$, from (\ref{K}) and (\ref{stirling}) it follows
$$
K_{n,[[y+ n]]}
\leq C \sqrt{\frac{n^{2n}}{y^y(2n-y)^{2n-y} }}=:\sqrt{F(n,y)}.
$$
A calculation shows that $n\mapsto F(n,y)$ is increasing on 
$I_{2}$ and that
\begin{equation}
	F_{1}(y):=\ln\left( \max_{n\in I_{2}}
	F(n,y)\right)=\ln F(|y|-|y|^{\frac{1}{2}+\delta})=-|y|^{2\delta}
	+O(|y|^{3\delta-1/2}).\nonumber 
\end{equation}
The case $n\in I_{3}$ is treated in a similar way. It follows that
\[\sum_{n\in I_{2}\cup I_{3}}K_{n,[\![y+ n]\!]}
\leq O(1)|y| e^{-C(\delta)|y|^{2\delta}}.\]
Finally let $n\in I_{4}$. Since the map $y\mapsto F(n,y)$ is increasing
when $|y|\leq n$, then $F(n,y)\leq F(n,n/2)\leq (2/3)^{n/2}$
for $n\in I_{4}$. Thus
$$
\sum_{n\in I_{4}}
K_{n,[\![y+n]\!]}=O(1)\sum_{n\geq 2|y|}\left(\frac{2}{3}\right)^{n/4}
=O(1)\left(\frac{2}{3}\right)^{|y|/2}.
$$
This completes the proof.
\end{proof}
In the following we will also need an analogous result for the
heat kernel $G$.
\begin{proposition}\label{tailsG}
Let  $y<0$, $\lambda>0$ and $\delta\in]0,1/2[$. Then, outside the interval
$I(|y|;\lambda,\delta)$ the integral of $G$ as well as any of its derivatives 
is transcendentally small, i.e.\ for any $k\geq 0$, we have
\begin{equation}\label{tailGesti}
\int_{\mathbb{R}^{+}\setminus I(|y|;\lambda,\delta)}
\left|\frac{\partial^k}{\partial x^k}G(t,y+\lambda t)\right|\,dt
\leq O(1)e^{-C(\delta,\lambda) |y|^{c(\lambda,\delta)}},  
\end{equation}
for some positive constants $C(\delta,\lambda)$ and $c(\lambda,\delta)$.
Moreover 
\begin{equation}\label{tailGesti1}
	\bigg|\sum_{n\notin I(|y|;\lambda,\delta)}
	\frac{\partial^k}{\partial x^k}G\left(n,
	y+\lambda n\right) \bigg|\leq
	O(1)e^{-C(\delta,\lambda) |y|^{c(\lambda,\delta)}},  
\end{equation}
\begin{equation}\label{tailGesti2}
	\bigg|\sum_{n\notin I(|y|;\lambda,\delta)}
	\frac{\partial^k}{\partial x^k}G\big(n,
	[\![y+\lambda n]\!]\big) \bigg|\leq
	O(1)e^{-C(\delta,\lambda) |y|^{c(\lambda,\delta)}},  
\end{equation}
\end{proposition}
\begin{proof}
 For simplicity, assume $\lambda=1$. Concerning the integral, we have
\begin{align*}
    \int_{0}^{|y|-|y|^{1/2+\delta}}
    \left|\frac{\partial^k}{\partial x^k}G(t,y+ t)\right|\,dt
    &=O(1)\int_{0}^{|y|-|y|^{1/2+\delta}}e^{-|y|^{2\delta}/5}\,dt
    \notag\\
    &=O(1)\cdot |y|\,e^{-|y|^{2\delta}/5},
\end{align*}
\begin{align*}
    \int_{|y|+|y|^{1/2+\delta}}^{\infty}
    &\left|\frac{\partial^k}{\partial x^k}G(t,y+ t)\right|\,dt
    =O(1)\cdot\left\{\int_{|y|^{1/2+\delta}}^{|y|}
    +\int_{|y|}^{\infty}\right\} e^{-\frac{\tau^{2}}{4(\tau+|y|)}}\,d\tau
    \notag\\
    &\leq C|y|e^{-|y|^{2\delta}/8}+\int_{|y|}^{\infty}e^{-\tau/8}\,d\tau
    =O(1)e^{-C|y|^{2\delta}}.
\end{align*}
Hence (\ref{tailGesti}) follows. The other two estimates can be obtained
from (\ref{tailGesti}), (\ref{2.13}) and (\ref{2.14}). 
\end{proof}
As a consequence of the previous propositions, we are authorized to add or
subtract the tails of the integrals/sums, introducing an error which is
exponentially decreasing with $y$. This will be frequently and tacitly
used in the following.
%
%
\subsection{Estimates on the approximation}
In order to estimate the variation of the second component of the DSPs we 
will need that $v^{(\lambda)}(x;\xi)$ is close, within
acceptable errors, to the corresponding function defined in terms of
the heat kernel. This function is given as
\[w^{(\lambda)}(x;\xi):=2\sum_{n=1}^\infty
G\left(\frac{n}{2},[\![z_n]\!]\right),\] 
where 
\begin{equation}\label{z_n}
	z_n=z_n(x,\lambda):=x+\lambda (n+1) -\xi.
\end{equation}
We will need to carefully keep track of
the dependence of $w^{(\lambda)}$ on $\lambda$.

We proceed estimate $w^{(\lambda)}$. Using  (\ref{2.10}),  (\ref{2.13}), 
(\ref{Gx_esti}) and (\ref{Gt_esti}),
for all $m\geq1$, $x\ll0$ and $\xi$ in the support of $\psi$, writing
$z=x-\xi+\lambda$, we have
\begin{align}
	w^{(\lambda)}(x;\xi)&=2\sum_{n=1}^\infty
	G\left(\frac{n}{2},[\![z_n]\!]\right)\nonumber\\
	&= 2\sum_{n=1}^\infty
	G\left(\frac{n}{2},z_n\right)- 2\sum_{n=1}^\infty
	\left\{G\left(\frac{n}{2},z_n\right) -
	G\left(\frac{n}{2},[\![z_n]\!]\right) \right\}\nonumber\\ 
	&= \frac{2}{\lambda}+O(1)|x|^{-m}+\nonumber\\ 
	&\quad - 2\sum_{n\in I(|z|;\lambda,\delta)}
	\left\{G_x\left(\frac{n}{2},z_n\right) \cdot
	(\!(z_n)\!) + O(1)n^{-3/2} 
	(\!(z_n)\!)^2\right\}\nonumber\\
	&= \frac{2}{\lambda} -2\sum_{n\in I(|z|;\lambda,\delta)}
	G_x\left(\frac{n}{2},z_n\right) \cdot
	(\!(z_n)\!) \ +\ O(1)|x|^{-1+\delta}, \label{general}
\end{align}
where $(\!(a)\!):=a-[\![a]\!]$ is the fractional part for any real number $a$.
In this calculation we have used that 
$| I(|z|;\lambda,\delta) |=O(1)|x|^{1/2+\delta}$ and that 
$n=O(1)|x|$ when $n\in I(|z|;\lambda,\delta)$. 
%
%
\paragraph{The case of rational speed} 
The estimate (\ref{general}) is valid for any speed 
$\lambda$. Now suppose that $\lambda$ is a rational speed,
\[\lambda=\frac{p}{q},\qquad p,q\in\mathbb N,\]
with $p$ and $q$ relatively prime, and consider  the sum on the right-hand side 
of (\ref{general}).
By writing $n=mq+j$ with $m\geq 0$, $j\in\{0,\dots,q-1\}$, 
Taylor expanding about the points $(t_m,x_m)$ where
$$t_m=mq/2, \qquad x_m=x_m(x,\lambda):=x+\lambda mq-\xi,$$
and 
using the formula
\begin{equation}\label{identity}
	\sum_{j=1}^{q}\left(\!\!\!\left(
	z+\frac{pj}{q}\right)\!\!\!\right)  
	=(\!(qz)\!)+\frac{q-1}{2},
\end{equation}
from (\ref{Gx_esti}) and (\ref{Gt_esti}), we obtain, for $z_n=z_n(x,\lambda)$ 
and $z=x-\xi+\lambda$,
\begin{align}
	&\sum_{n\in I(|z|;\lambda,\delta)}
	G_x\left(\frac{n}{2},z_n\right) (\!(z_n)\!) \nonumber\\
	&= \underbrace{\sum_{m\geq
	0}\sum_{j=0}^{q-1}}_{mq+j \in I(|z|;\lambda,\delta)} \nquad
	G_x\left(\frac{mq+j}{2},x+\lambda(mq+j+1)-\xi\right)
	(\!(z_{mq+j})\!)  \nonumber\\
	&= \underbrace{\sum_{m\geq
	0}\sum_{j=0}^{q-1}}_{mq+j \in I(|z|;\lambda,\delta)} \nquad
	\left\{ G_x\left(\frac{mq}{2},x_m\right) 
	+\sup_{t,x}\big(|G_{xt}|+|G_{xx}|\big)O(q)\right\}  (\!(z_{mq+j})\!)  \nonumber\\
	&= \sum_{m\geq 1} G_x\left(\frac{mq}{2},x_m\right)
	\sum_{j=0}^{q-1}(\!(z_{mq+j})\!) \nonumber\\ 
	&\qquad\qquad\qquad\qquad+
	O(q)|x|^{1/2+\delta}\big(|x|^{-2} + 
	|x|^{-3/2}\big)+ e^{-C|x|}  \label{blabla} \\[4pt]
	&= \left\{(\!(q(x-\xi) )\!) + \frac{q-1}{2}\right\}
	\sum_{m\geq 1} 
	G_x\left(\frac{mq}{2},x_m\right) + O(q)|x|^{-1+\delta}.\nonumber
\end{align}
In this calculation we have used that $| I(|z|;\lambda,\delta) |=O(1)|x|^{1/2+\delta}$ and that 
$n=O(1)|x|$ when $n\in I(|z|;\lambda,\delta)$.  
Analogously to (\ref{2.14}), one can prove that
\[
\sum_{m\geq 1} G_x\left(\frac{mq}{2}, x_m\right)=O(q^{M-1}) |x|^{-\frac{M+1}{2}+\delta},
\]
for every integer $M$. We thus have that
\begin{align}
	\sum_{n\in I(|z|;\lambda,\delta)}
	&G_x\left(\frac{n}{2},z_n\right) \cdot (\!(z_n)\!) \notag\\
	=& \left\{(\!(q(x-\xi))\!) + \frac{q-1}{2}\right\}\sum_{m\geq 1} 
	G_x\left(\frac{mq}{2}, x_m\right) + O(q)|x|^{-1+\delta}\notag\\
	=& O(q^{M}) |x|^{-\frac{M+1}{2}+\delta}+O(q) |x|^{-1+\delta}. \label{w_esti2}
\end{align}
From (\ref{general}) and (\ref{w_esti2}) we conclude that when the speed 
$\lambda=p/q$ is a rational we have 
\begin{equation}\label{w_lambda}
	w^{(\lambda)}(x;\xi)=\frac{2}{\lambda} + O(q^{M}) |x|^{-\frac{M+1}{2}+\delta}+O(q) |x|^{-1+\delta}.
\end{equation}
Due to the dependence on the denominator $q$, equation (\ref{w_lambda}) is not useful for computing
the variation of differences $w^{(\lambda)}-w^{(\tilde\lambda)}$ as $\tilde\lambda\to\lambda$.
Instead, (\ref{w_lambda}) will be used for a fixed reference speed 
$\lambda$, while we need  an alternative analysis to estimate $w^{(\tilde{\lambda})}$, where 
$\tilde{\lambda}=\lambda+\eps$ is a small perturbation of $\lambda$.

For this we return to the right hand side of (\ref{blabla}), with $x_m=x_m(x,\tilde{\lambda})$
and $z_{n}=z_{n}(x,\tilde{\lambda})$.
We establish the following technical result. 
\begin{proposition}\label{technical}
	Let $\lambda=p/q$ and $\tilde{\lambda}=\lambda+\eps$, with $|\eps|\ll 1$. Then
	for $z=x-\xi$, where $\xi$ lies in the support of $\psi$, there holds
	\begin{align}
	\sum_{m\geq 0}
	&G_x\left(\frac{mq}{2},z+\tilde\lambda mq\right) \cdot 
	\sum_{j=1}^{q} (\!(z+mq\eps+\tilde\lambda j)\!) \nonumber \\
	&=\int_0^\infty G_x\left(\frac{sq}{2},z+\tilde\lambda sq
	\right) \cdot (\!(q(z+sq\eps))\!)\, ds + \nonumber \\
	&\qquad\qquad\qquad\qquad\qquad +O(q)|x|^{-1+\delta}
	+O(|\eps| q^2)|x|^{-1/2+\delta}  \label{tech}
	\end{align}
\end{proposition}
\begin{proof}
We compare both the sum and the integral to the same integral and
show that in each case the error is 
$O(q)|x|^{-1+\delta}+O(|\eps| q^2) |x|^{-1/2+\delta}$. Let 
$a_j=z+\tilde\lambda j$ and consider the following difference 
\begin{align} 
	&\left|\sum_{m\geq 0} G_x\left(\frac{mq}{2},z+\tilde\lambda
	mq\right) \cdot \sum_{j=1}^{q} (\!(z+mq\eps+\tilde\lambda j)\!)
	\right.- \nonumber \\  
	&\left. \qquad\qquad\qquad-\int_0^\infty
	G_x\left(\frac{sq}{2},z+\tilde\lambda sq  
	\right) \cdot \sum_{j=1}^{q} (\!(a_j+sq\eps)\!) \,
	ds\right| \nonumber \\ 
	&\leq \left|\sum_{m\geq 0} \int_m^{m+1}
	\left\{G_x\left(\frac{mq}{2},z+\tilde\lambda mq\right)
	\right.\right.- \nonumber \\
	&\left.\left.\qquad\qquad\qquad
	-G_x\left(\frac{sq}{2},z+\tilde\lambda sq 
	\right)\right\} \sum_{j=1}^{q} (\!(a_j+mq\eps)\!)  \, ds\right| +\nonumber \\  
	&+\left|\sum_{m\geq 0} \int_m^{m+1}
	G_x\left(\frac{sq}{2},z+\tilde\lambda sq \right) \left\{
	\sum_{j=1}^{q} (\!(a_j+mq\eps)\!)-
	(\!(a_j+sq\eps)\!) \right\}\,
	ds\right| \nonumber \\ 
	&\leq \underbrace{\sum_{m\geq 0}\sum_{j=1}^{q}}_{mq+j 
	\in I(|z|;\tilde\lambda,\delta)}
	\left\{\sup_{t,x}\big(|G_{xt}|+|G_{xx}|\big)O(q)\right\} 
	(\!(a_j+mq\eps)\!)  +e^{-C|x|}+  \nonumber \\ 
	&\quad+\frac{O(1)}{|x|} \sum_{j=1}^{q}\sum_{mq\in I(|z|;\tilde\lambda,\delta)}
	 \int_m^{m+1} \big|(\!(a_j+mq\eps)\!)-(\!(a_j+sq\eps)\!)\big|\,ds. 
	 \label{blabla2}
\end{align}	
The first sum on the right hand side of (\ref{blabla2}) is bounded by 
$O(q)|x|^{-1+\delta}$. To estimate the second sum, we divide the interval 
$I(|z|;\tilde\lambda,\delta)$ into 
$r=O\big(|\eps| q \cdot| I(|z|;\tilde\lambda,\delta)| \big)=O\big(
|\eps| q  \big) |x|^{1/2+\delta}$ intervals $J_1,\dots, J_r$ of equal length 
$ 1/(|\eps| q)$, and re-write the sum over $m$ as
\[
\sum_{k=1}^r \int_{J_k} \big|(\!(a_j+[\![s]\!]q\eps)\!)-(\!(a_j+sq\eps)\!)\big|\,ds.
\]
The integrand is bounded by $|\eps| q$ except on a sub-interval of length $1$ 
where it is $O(1)$. 
Hence, for each $k$, the integral over $J_k$ is bounded by order one. It follows
that the second sum on the right hand side of (\ref{blabla2}) is bounded by 
$O(|\eps| q^2) |x|^{-1/2+\delta}$.

On the other hand, using (\ref{identity}), (\ref{2.7}) and (\ref{Gx_esti}), 
we have 
\begin{align}
	&\left|\int_0^\infty
	G_x\left(\frac{sq}{2},z+\tilde\lambda sq 
	\right) \cdot (\!(q(z+sq\eps))\!)\, ds \right.- \notag\\
	&\qquad\qquad-\left.\int_0^\infty 
	G_x\left(\frac{sq}{2},z+\tilde\lambda sq 
	\right) \cdot \sum_{j=1}^{q} (\!(a_j+sq\eps)\!) \, ds \right|\notag\\
	&= \left|\int_0^\infty G_x\left(\frac{sq}{2},z+\tilde\lambda sq
	\right) \cdot \left\{(\!(q(z+sq\eps)\!)) -\sum_{j=1}^{q}
	(\!(a_j+sq\eps)\!)\right\}\, ds \right|\notag\\ 
	&=\left| \int_0^\infty G_x\left(\frac{sq}{2},z+\tilde\lambda sq
	\right)\cdot
	\phantom{\left\{ \sum_{j=1}^q \right\}}
	\right.\notag\\
	&\qquad\qquad\cdot \left\{ \sum_{j=1}^q \Big( (\!(z+\lambda
	j+sq\eps)\!) -(\!(a_j+sq\eps)\!)\Big)  -\frac{q-1}{2} \right\}\,
	ds\Bigg|\notag\\  
	&= \left|\int_0^\infty 
	G_x\left(\frac{sq}{2},z+\tilde\lambda sq \right) \cdot 
	\sum_{j=1}^{q} \Big( (\!(z+\lambda j+sq\eps)\!)-(\!(a_j+sq\eps)\!)
	\Big) \, ds\right| \notag\\   
	&\leq \frac{O(1)}{q|z|}\sum_{j=1}^{q} \int_{I(|z|;2\tilde\lambda,\delta)}
	\Big| (\!(z+\lambda j+2\eps t)\!)-(\!(a_j+2\eps t)\!) \Big|\, dt
	+O(1)|x|^{-M} \notag\\
	&=  \frac{O(1)}{ q|\eps z|}\sum_{j=1}^{q} \int_{J(|z|;j,\eps,\tilde\lambda,\delta)}
	\big|(\!(\eta)\!)-(\!(\eta+\eps j)\!)\big|\, d\eta+O(1)|x|^{-M},
	\label{blabla3}
\end{align}
where $J(|z|;j,\eps,\tilde\lambda,\delta)=z+\lambda j +2\eps I(|z|;2\tilde\lambda,\delta)$.
Each of the integrals in the sum is therefore of order $O(\eps^2 q)|z|^{1/2+\delta}$. Finally,
using that $z=x-\xi$ and that $\xi$ varies in a compact interval, 
from (\ref{blabla2}) and (\ref{blabla3}) we thus obtain (\ref{tech}).
\end{proof}
Using (\ref{blabla}) and (\ref{general}) for speed $\tilde\lambda$ together with 
Proposition \ref{technical} yields
\begin{align*}
	 w^{(\tilde\lambda)}(x;\xi) &= \frac{2}{\tilde\lambda} -2\sum_{m\geq 0}
	G_x\left(\frac{mq}{2},z+\tilde\lambda mq\right) \cdot 
	\sum_{j=1}^{q} (\!(z+mq\eps+\tilde\lambda j)\!)+ \\
	&\qquad\qquad\qquad+O(q)|x|^{-1+\delta} \\[4pt]
	&=\frac{2}{\tilde\lambda} - 2\int_0^\infty
	G_x\left(\frac{sq}{2},z+\tilde\lambda sq\right) \cdot
	(\!(q(z+sq\eps))\!)\, ds\\ 
	&\qquad\qquad\qquad+O(q)|x|^{-1+\delta}+O(|\eps| q^2)|x|^{-1/2+\delta} 
	\\[4pt]
	&=\frac{2}{\tilde\lambda} -
	\frac{2}{q} \int_0^\infty G_x\left(\frac{t}{2},z+\tilde\lambda t\right) 
	\cdot  (\!(qz+q\eps t)\!)\, dt+E(|x|;\eps,q),
\end{align*}
where $z=x-\xi$ and 
\begin{equation}\label{error}
	E(y;\eps,q)=\frac{O(q)}{y^{1-\delta}}+\frac{O(|\eps| q^2)}{y^{1/2-\delta}}.
\end{equation}
We proceed by simplifying this expression. Put $y=-x\gg0$, such that
$z+\tilde \lambda t=-y-\xi+\tilde \lambda t$, and introduce the
coordinate $\tau$ by   
\[\tau\sqrt{\frac{y}{\tilde\lambda}}=y-\tilde\lambda t.\] 
Thus $\tau$-values outside an interval $[-Cy^\delta,Cy^\delta]$
corresponds to $t$-values for which the contribution to the integral
is exponentially small. A simple calculation now shows that
\[G_x\left(\frac{t}{2},z+\tilde \lambda t\right) =
\frac{1}{\sqrt{2\pi}} \frac{\tilde\lambda}{y}\tau e^{-\tau^2/2} +
O(1) y^{-3/2+\delta}.\] 
Substitution into the expression above yields
\begin{align}
w^{(\tilde\lambda)}(&-y;\xi)=\frac{2}{\tilde\lambda}
+E(y;\eps,q)-  \nonumber \\
& -\frac{1}{q}\sqrt{\frac{2}{\pi\tilde\lambda y}}
{\underset{|\tau|<Cy^\delta}{\int}} \tau e^{-\tau^2/2}
\cdot  \left(\!\!\! \left(\frac{q\eps}{\tilde\lambda}
\left( y-\tau\sqrt{\frac{y}{\tilde\lambda}} \right) -qy-q\xi 
 \right)\!\!\! \right)\, d\tau.\label{w_lambda_tilde}
\end{align} 
Recalling (\ref{w_lambda}) we get that for a rational speed $\lambda=p/q$ 
and a perturbation $\tilde\lambda=\lambda+\eps$, there holds
\begin{align}
	w^{(\lambda)}(x;&\xi)-w^{(\tilde\lambda)}(x;\xi) =
	\frac{2}{\lambda} -\frac{2}{\tilde\lambda} + E(y,\eps,q)-\nonumber \\ 
	& -\frac{C_1}{qy^{1/2}}
	{\underset{|\tau|<Cy^\delta}{\int}} \tau e^{-\tau^2/2} \cdot
	\left(\!\!\! \left(q\xi+qy-\frac{q\eps}{\tilde\lambda} 
	\left( y-\tau\sqrt{\frac{y}{\tilde\lambda}} \right) 
	\right)\!\!\! \right) \, d\tau, \label{w_est}
\end{align}
where $C_1$ is a positive constant.  We will use this estimate for
rational speeds $\tilde\lambda$ converging to a fixed rational speed $\lambda$.
 The relevance of the estimate is that it does not depend explicitly on the 
 denominator of the approximating speeds.
%
%
\subsection{Estimating $v^{(\lambda)}-w^{(\lambda)}$} 
In this subsection we show that $w^{(\lambda)}$ is a sufficiently good approximation of 
$v^{(\lambda)}$. We have the following estimate which actually holds for any real speed $\lambda$.
\begin{proposition}\label{v^lambda_approx}
	Provided $\xi$ is in the compact support of $\psi$, 
	the functions $v^{(\lambda)}(x;\xi)$ and
	$w^{(\lambda)}(x;\xi)$ satisfy
	\begin{equation}\label{v_approx}
		v^{(\lambda)}(x;\xi) - w^{(\lambda)}(x;\xi) = 
		\sum_{\underset{n-[\![z_n]\!] \mathrm{odd}}{n\in I(|z|;\lambda,\delta)}}
	2G_x\left(\frac{n}{2}, z_n \right) + O(1)|x|^{-1+2\delta},
	\end{equation}
	when $x<0$, where $z=x-\xi+\lambda$ and $z_n=z_n(x,\lambda)$ 
	is given by (\ref{z_n}).
\end{proposition}
\begin{proof}
Observe that
\[G\left(\frac{n}{2}, [\![z_n]\!]+1\right)=G\left(\frac{n}{2}, [\![z_n]\!]\right)+ 
G_x\left(\frac{n}{2}, z_n\right) + O(1) n^{-3/2}.\]
By Proposition \ref{tails} only indices $n\in
I(|z|;\lambda,\delta)$ in the  
sums defining $v^{(\lambda)}$ and $w^{(\lambda)}$ are significant.
If $n\in I(|z|;\lambda,\delta)$ then $n=O(1)|x|$ while $[\![z_n]\!]$ is $O(1)|x|^{1/2+\delta}$, so that we are within the
validity of the approximation (\ref{approx}). 
Since exactly one of $n-[\![z_n]\!]$, $n-([\![z_n]\!]+1)$ is even, we have
\begin{align*}
	v^{(\lambda)}(x;\xi) &= \sum_{n\in I(|z|;\lambda,\delta)} \big(
	K_{n,[\![z_n]\!]}+ K_{n,[\![z_n]\!]+1}\big)+O(1)e^{-C|x|} \\
	&= \sum_{\underset{n-[\![z_n]\!] \mathrm{even}}{n\in
	I(|z|;\lambda,\delta)}} 2G\left(\frac{n}{2}, [\![z_n]\!]\right) 
	+\sum_{\underset{n-[\![z_n]\!] \mathrm{odd}}{n\in I(|z|;\lambda,\delta)}}
	2G\left(\frac{n}{2}, [\![z_n]\!]+1\right)+ \\
	&\quad+\sum_{n\in I(|z|;\lambda,\delta)} O(1)n^{-3/2+\delta}
	+O(1)e^{-C|x|} \\ 
	&= \sum_{n\in
	I(|z|;\lambda,\delta)} 2G\left(\frac{n}{2}, [\![z_n]\!]\right) 
	+\sum_{\underset{n-[\![z_n]\!] \mathrm{odd}}{n\in I(|z|;\lambda,\delta)}}
	2G_x\left(\frac{n}{2}, z_n\right)+ \\
	&\quad+\sum_{n\in I(|z|;\lambda,\delta)} O(1)n^{-3/2+\delta}
	+O(1)e^{-C|x|} \\ 
	&=w^{(\lambda)}(x;\xi)+ \sum_{\underset{n-[\![z_n]\!]
	\mathrm{odd}}{n\in I(|z|;\lambda,\delta)}}
	2G_x\left(\frac{n}{2}, z_n\right)  +O(1)|x|^{-1+2\delta},
\end{align*} 
since $I(|z|;\lambda,\delta)$ has length $O(1)|x|^{1/2+\delta}$. 
\end{proof}
To estimate the sum on the right-hand side of (\ref{v_approx}) we will make 
particular choices for the speeds $\lambda$ and $\tilde\lambda$. In what follows we will fix 
\begin{equation}\label{speeds}
	\lambda=\frac{1}{2},\qquad\tilde\lambda=\frac{k}{2k+1},
\end{equation}
such that $\eps=-\textstyle\frac{1}{4k+2}$. In order to satisfy the non-resonance condition of 
Majda and Ralston \cite{MajdaRalston79} we will let $k$ be an even integer.
We have the following key estimates that are independent of $k$ (and $\eps$).
\begin{proposition}\label{key_sum}
	With $\lambda$ and $\tilde\lambda$ given by (\ref{speeds}) we have
	 \begin{equation}\label{key1}
		\sum_{\underset{n-[\![z_n(x,\lambda)]\!] \mathrm{odd}}{n\in I(|z|;\lambda,\delta)}}
		G_x\left(\frac{n}{2}, z_n(x,\lambda) \right) = \frac{O(1)}{|x|^{1-\delta}},
	\end{equation}
	and 
	 \begin{equation}\label{key2}
		\sum_{\underset{n-[\![z_n(x,\tilde\lambda)]\!] \mathrm{odd}}{n\in I(|z|;\tilde\lambda,\delta)}}
		G_x\left(\frac{n}{2}, z_n(x,\tilde\lambda) \right) = \frac{O(1)}{|x|^{1-\delta}}.
	\end{equation}
\end{proposition}
\begin{proof}
We give the proof for (\ref{key2}), the proof of (\ref{key1}) being similar.
Given $\xi$ and $x$, we set $z=x-\xi+\tilde\lambda$ such that $z_n(x,\tilde\lambda)=z+\tilde\lambda n$. 
We denote by $\mathcal N$ the set of the integers $n$ contributing to the sum (\ref{key2}). 
Let a $k$-block denote a half-open 
interval of length $2k+1$ consisting of $k$ consecutive subintervals of equal length 
$1/\tilde\lambda=2+1/k$, each on which the function $B(s):= [\![z+\tilde\lambda s]\!]$ is constant. 
Without loss of generality we can assume that 
$I(|z|;\tilde\lambda,\delta)$ is exactly partitioned into finitely many 
$k$-blocks.
A $k$-block contains $2k+1$ integers and each subinterval has length $2+1/k$. It follows that in each
$k$-block there are $k-1$ subintervals containing exactly two integers, and one subinterval containing
three integers.  Since $B(s)$ is constant on each of the subintervals, the function
$n\mapsto n-  [\![z+\tilde\lambda n]\!]$ takes on one even and one odd value on each subinterval that 
contains exactly two integers. As $k$ is even two consecutive $k$-blocks will
contain exactly $2k+1$ integers $n$ in $\mathcal N$. Thus, of two 
consecutive $k$-blocks, one contains $k$, and the other $k+1$, of integers $n$ in $\mathcal N$. 
We observe that in the sequence of these indices $n$, 
the elements are at most a distance $3$ apart from each  other. Finally, 
since $k$ is an even integer, in all subsequent unions of two consecutive $k$-blocks, 
the distribution of the indices $n\in \mathcal N$ is the same. 
We can therefore define a map $\mu$ that maps $\mathcal N$ bijectively onto the {\em regular}
grid of even integers in $I(|z|;\tilde\lambda,\delta)$, in such a way that $|n-\mu(n)|\leq 3$.
The sum in (\ref{key2}) can therefore be estimated as follows,
\begin{align*}
    \sum_{n\in\mathcal N} G_x\left(\frac{n}{2}, 
    \right.&z_n(x,\tilde\lambda) \bigg)
    =\sum_{n\in\mathcal N} G_x\left(\frac{\mu(n)}{2}, z_{\mu(n)}(x,\tilde\lambda)
    \right)+\\
    &\qquad+\sum_{n\in\mathcal N} \left[G_x\left(\frac{n}{2}, z_n(x,\tilde\lambda) \right)
    -G_x\left(\frac{\mu(n)}{2}, z_{\mu(n)}(x,\tilde\lambda) \right)\right]\\
    &=\sum_{n\in 2\mathbb Z\cap I(|z|;\tilde\lambda,\delta)} 
    G_x\left(\frac{n}{2}, z_n(x,\tilde\lambda) \right) +\\
    &\qquad\qquad+ O(1)\sup(|G_{xt}|+|G_{xx}|) \cdot|I(|z|;\tilde\lambda,\delta)| \\
    &= \frac{O(1)}{|x|^{1-\delta}},
\end{align*}
where we have used (\ref{2.7}), (\ref{Gx_esti}) and (\ref{Gt_esti}).
\end{proof}
\begin{remark}
	We note that the complexity of the preceding arguments is essentially due to the fact that we
	are working with the Lax-Friedrichs scheme. The same computations would be significantly
	simpler for e.g. the upwind scheme, in which case the complication of even and odd terms do not 
	occur.
\end{remark}
Combining Proposition \ref{v^lambda_approx} and Proposition \ref{key_sum} we conclude that, if 
the velocities $\lambda$ and $\tilde\lambda$ are given by (\ref{speeds}), then
 \begin{equation}\label{key_est}
	v^{(\lambda)}(x;\xi) - v^{(\tilde\lambda)}(x;\xi)= 
	w^{(\lambda)}(x;\xi)  - w^{(\tilde\lambda)}(x;\xi) + \frac{O(1)}{|x|^{1-\delta}}.
\end{equation}
%
%
\section{Variation of $V^{(\lambda)}-V^{(\tilde\lambda)}$} 
In the remaining part of the paper the velocities $\lambda$ and $\tilde\lambda=\lambda +\eps$ 
are given by (\ref{speeds}).
Recalling that $\psi$ has  compact support we conclude from 
(\ref{repr}), (\ref{key_est}) and (\ref{w_est}) that 
\begin{align}
	&V^{(\lambda)}(x)-V^{(\tilde\lambda)}(x)= 
	-\frac{\eps}{\lambda\tilde\lambda}\int_{-\infty}^\infty\psi(\xi)\,
	d\xi +E(y,\eps)+  \nonumber  \\	
	&\ +\frac{C_1}{y^{1/2}}
	{\overset{\quad\infty}{\underset{-\infty}{\int}}}  
	{\underset{|\tau|<Cy^\delta}{\int}} \psi(\xi)\tau 
	e^{-\tau^2/2} \left(\!\!\! \left(2\xi+2y-\frac{2\eps}{\tilde\lambda} 
	\left( y-\tau\sqrt{\frac{y}{\tilde\lambda}} \right) 
	\right)\!\!\! \right)\, d\tau\, d\xi ,\nonumber\\ 
	&= A(\eps)+E(y,\eps)+\frac{C_1}{y^{1/2}}
	\int_{|\tau|<Cy^\delta} \tau
	e^{-\tau^2/2}h(\tau;y,\eps) \, d\tau,  \label{Delta_V}
\end{align}
where $y=-x$ and where we have introduced 
\begin{align*}
	E(y,\eps) &:= E(y,\eps,2)\int_{-\infty}^\infty\psi(\xi)\,  d\xi,\\[4pt]
	A(\eps)&:=-\frac{\eps}{\lambda\tilde\lambda}
	\int_{-\infty}^\infty\psi(\xi)\,  d\xi, \\[4pt]
	h(\tau;y,\eps) &:= \int_{-\infty}^\infty \psi(\xi)
	\left(\!\!\! \left(2\xi+2y-\frac{2\eps}{\tilde\lambda} 
	\left( y-\tau\sqrt{\frac{y}{\tilde\lambda}} \right) 
	\right)\!\!\! \right)  \, d\xi.  
\end{align*}
Making the change of variables $\eta=2\xi$ and denoting 
$\phi(\eta)=\psi(\eta/2)$, $\beta=2\lambda
/\tilde\lambda$, and $\gamma=2/\tilde\lambda^{3/2}$, we 
have that
\begin{equation}\label{h}
	h(\tau;y,\eps)=\frac{1}{2}\int_{-\infty}^\infty \phi(\eta)
	\big(\!\big(\eta+\beta y + \gamma\eps\sqrt{y}\tau\big)\!\big)\, d\eta.
\end{equation}
We proceed to use (\ref{Delta_V}) to  show that the function
$V^{(\lambda)}-V^{(\tilde\lambda)}$ contains an $O(1)$ amount of
variation on an interval of the form
\[J(\eps):=[-C\eps^{-2/(1+2\delta)},-c\eps^{-2/(1+2\delta)}].\]
Recalling (\ref{error}) we see that 
$E(y,\eps)$ is of order $O(\eps^{2(1-\delta)/(1+2\delta)})$ on $J(\eps)$. 

As a first step we consider the limiting case where $\phi(\eta)/2$ is the Dirac delta function 
centered at $\eta=0$. In this case 
\[h(\tau;y,\eps)=h_0(\tau;y,\eps):=\big(\!\big(\beta y + \gamma\eps\sqrt{y}\tau\big)\!\big),\]
such that 
\begin{equation}\label{expr}
	V^{(\lambda)}(x)-V^{(\tilde\lambda)}(x)= A(\eps) 
	+ E(y,\eps) +\frac{C_1}{y^{1/2}}H_0(y;\eps),
\end{equation}
where 
\begin{equation}\label{H}
    H_0(y;\eps):=\int_{|\tau|<Cy^\delta} \tau 
    e^{-\tau^2/2}\big(\!\big(\beta y + \gamma\eps\sqrt{y}\tau\big)\!\big)\, 
    d\tau.
\end{equation}
\begin{lemma}\label{lemma}
	There exist an $O(\eps^{-1/(1+2\delta)})$ number of points $y_1,\dots,y_L$ in $-J(\eps)$
	with
	\begin{equation}
		H_0(y_n;\eps)=\left\{\begin{array}{ll}
		-O(1) & \mbox{for $n$ odd,}\\[2pt]
		O(\eps^{2\delta/(1+2\delta)}) & \mbox{for $n$ even.}\end{array}\right.
	\end{equation}
\end{lemma}
\begin{proof}
    Set 
    \[y_1:= [\![\eps^{-2/(1+2\delta)}]\!]/\beta,\]
    and let
    \[y_n := y_1 + \frac{1}{\beta}\left(n[\![\eps^{-1/(1+2\delta)}]\!]+\frac{n+1}{2}\right),\]
    for $n=2,\dots,[\![\eps^{-1/(1+2\delta)}]\!]$. These choices imply that
    \[\left.\big(\!\big(\beta y_n + \gamma\eps\sqrt{y_n}\tau\big)\!\big)\right|_{\tau =0}
    = (\!(\beta y_n)\!)=\left\{\begin{array}{ll}
    0 & \mbox{for $n$ odd,}\\
    1/2 & \mbox{for $n$ even.}\end{array}\right.\] 
    We observe that an $O(1)$ change in the  constant $C$ in (\ref{H}) induces an exponentially small 
    error (with respect to $y$). In order to simplify the computations we make the following choices
    for the constant $C$ in (\ref{H}),
    \[C=C_n:=\left\{\begin{array}{ll}
	\frac{1}{\gamma\eps y_n^{\delta+1/2}}=O(1) & \mbox{for $n$ odd,}\\[6pt]
	\frac{1}{2\gamma\eps y_n^{\delta+1/2}}=O(1) & \mbox{for $n$ even.}\end{array}\right.\] 
	Thus, up to an exponentially small error, we have for $n$ odd
    \begin{align}\label{H_est}
    H_0(y_n;\eps) &= \int_{|\tau|<C_n y_n^\delta} \tau
    e^{-\tau^2/2}\big(\!\big(\beta y_n + \gamma\eps\sqrt{y_n}\tau\big)\!\big)\, 
    d\tau \notag \\
    &=\int_{-\infty}^0\tau e^{-\tau^2/2}=-O(1),
    \end{align}
    while for $n$ even a similar argument yields
    \[H_0(y_n;\eps) = O(1)\eps^{2\delta/(1+2\delta)}.\]
\end{proof}
For a fixed, small $\eps$ it now follows from Lemma \ref{lemma} and (\ref{expr}) that the function 
$V^{(\lambda)}(x)-V^{(\tilde\lambda)}(x)$ alternates between values 
\[ A(\eps) +O(1)\eps^{2(1-\delta)/(1+2\delta)}-O(1)\eps^{1/(1+2\delta)},\]
for $x=-y_n,$ $n$ odd, and 
\[A(\eps) +O(1)\eps^{2(1-\delta)/(1+2\delta)}+O(1)\eps,\]
for $x=-y_n,$ $n$ even. 
Provided $\eps$ is small enough, this implies that in every interval 
$[-y_{n+1},-y_{n}]$, the function $V^{(\lambda)}(x)-V^{(\tilde\lambda)}(x)$ 
contains an $O(1)\eps^{1/(1+2\delta)}$ amount of variation. Since there are
an $O(1)\eps^{-1/(1+2\delta)}$ number of such intervals,
this argument shows that, in the limiting case where $\phi(\eta)/2$ in
(\ref{h}) is a Dirac delta, the function $V^{(\lambda)}(x)-V^{(\tilde\lambda)}(x)$
contains at least an $O(1)$ amount of variation on $J(\eps)$.

It remains to argue that the same result holds whenever $\phi(\eta)/2$ is close 
to a Dirac delta function. For this it is sufficient to show that the difference 
$|H(y;\eps)-H_0(y;\eps)|$ can be made arbitrarily small independently of the 
$O(1)$-estimate in (\ref{H_est}). This will be accomplished by imposing that the
support of $\phi$ (hence $\psi$) be sufficiently small.
We have 
\begin{align*}
    &|H(y;\eps)-H_0(y;\eps)| \\
    &\leq \int_{-\infty}^{+\infty}\frac{\phi(\eta)}{2} \int_{|\tau|<C y^\delta}
    \left|K(\tau) \Big\{\big(\!\big(F(y;\eps,\tau)-\eta\big)\!\big)-
    \big(\!\big(F(y;\eps,\tau)\big)\!\big)\Big\}\right|\, d\tau\,d\eta,
\end{align*}
where $F(y;\eps,\tau)=\beta y + \gamma\eps\sqrt{y}\tau$. Since $\eps\sqrt{y}$ is $O(1)$ 
it follows that there is a constant $\bar C$ such that
\[\int_{|\tau|<C y^\delta} \Big|K(\tau) \Big\{
	\big(\!\big(F(y;\eps,\tau)-\eta\big)\!\big)-\big(\!\big(F(y;\eps,\tau)
	\big)\!\big)\Big\} \Big|\, d\tau \leq \bar C\eta.\]
We thus have that
\[|H(y;\eps)-H_0(y;\eps)|\leq C\int_{-\infty}^{+\infty}\phi(\eta)\eta\, d\eta,\]
which can be made arbitrarily small, compared to the $O(1)$-estimate in
(\ref{H_est}), by choosing $\phi$ sufficiently close to a Dirac delta.

We have thus proved the following theorem.
\begin{theorem}
	Discrete shock profiles for the Lax-Friedrichs scheme for strictly 
	hyperbolic systems of two conservation laws of the form (\ref{eq1})-(\ref{eq2}) do not depend 
	continuously in BV on their speeds. More precisely, one can find a sequence of rational 
	speeds $\lambda_n$ converging to $\lambda\in\mathbb Q$, for which there are  discrete 
	shock profiles $\Psi_n$ and $\Psi$ of speeds $\lambda_n$ and $\lambda$, 
	respectively, and such that 
	\[\mbox{T.V.} (\Psi_n-\Psi)=O(1),\]
	independently of $n$.
\end{theorem}
\begin{remark}
	From (\ref{repr}), (\ref{v_approx}), (\ref{key1}) and (\ref{w_lambda}),
	we observe that with our particular choices of speeds we have 
	\[V^{(\lambda)}(x+\Delta x) - V^{(\lambda)}(x)
	= \frac{O(1)}{|x|^{1-\delta}},\]
	for $x\ll 0$ and $\Delta x=O(1)$. It follows that an $O(1)$ translation of $V^{(\lambda)}$ 
	relative to $V^{(\tilde\lambda)}$ changes the total variation of their 
	difference only by $O(\eps)$ in the region $J(\eps)$. 
\end{remark}

\end{document}